\documentclass[11pt]{article}
\title{Harmonic Splittings of Surfaces}
\author{Benson Farb\thanks{Partially supported by NSF grant DMS 9704640} 
\and Michael Wolf\thanks{Partially 
supported by NSF grants
DMS 9300001 and DMS 9626565}}
\date{Final Version, June 8, 1999}

\newtheorem{theorem}{Theorem}[section]
\newtheorem{proposition}[theorem]{Proposition}
\newtheorem{lemma}[theorem]{Lemma}

\def\proof{{\bf {\medskip}{\noindent}Proof: }}

\def\remark{{\bf {\bigskip}{\noindent}Remark. }}

\def\endproof{$\diamond$ \bigskip}

\def\title{\em}

\def\bar{\overline}

\def\ree{\mbox{Re}}

\def\qd{\mbox{QD}}

\def\wt{\widetilde}
\def\tl{\tilde}
\def\tec{Teichm\"uller\ }

\newcommand\z{\zeta}
\newcommand\vp{\varphi}
\newcommand\om{\omega}
\newcommand\cd{\cdot}
\newcommand\SF{\mathcal{F}}

\newcommand\SN{\mathcal{N}}
\newcommand\SR{S}

\newcommand\hyp{\mbox{\bf H}}
\newcommand\C{\mbox{\bf C}}

\newcommand\R{\mbox{\bf R}}

\newcommand\Q{\mbox{\bf Q}}

\renewcommand\to{\rightarrow}
\newcommand\x{\times}
\newcommand{\ox}{\ensuremath{\otimes}}
\newcommand\MF{{\cal M}{\cal F}}
\newcommand\PMF{{\cal P}{\cal M}{\cal F}}
\newcommand\Teich{\mbox{Teich}}
\newcommand\ov{\overline}
\newcommand\p{\partial}
\newcommand \pov {\partial \overline}
\begin{document}

\maketitle

\section{Introduction}

Let $\Gamma=\pi_1(S)$ be the fundamental group of a closed surface $S$ of 
genus at least two. Morgan-Shalen showed \cite{MS2}, \cite{GiSh} 
that every point in the Thurston compactification $\PMF(S)$ of the
Teichmuller space $\Teich(S)$ gives an isometric  $\Gamma$-action on an 
$\R$-tree.  Given a measured foliation ${\cal F}\in \PMF(S)$, the action is
simply the $\Gamma$-action on the leaf space  of the lift of $\cal F$ to
$\hyp^2$.  This action is  {\em small} in the sense that edge stabilizers 
do not contain rank two free groups.   It is also {\em minimal} in the
sense that it leaves no proper subtree invariant.

Shalen \cite{Sh} conjectured that every minimal 
small action of $\Gamma$ on an
$\R$-tree $T$ arises in this way.  This problem has several important
applications in low-dimensional geometry and topology (see \cite{Ot}).  
Partial results were obtained by Morgan-Shalen \cite{MS1} and
Gillet-Shalen \cite{GiSh}.  

The conjecture was eventually proven in two steps: Morgan-Otal \cite{MO} 
(see also \cite{Ha}) constructed the candidate foliation, 
with
dual $\R$-tree $R$, and a $\Gamma$-equivariant 
morphism $\phi :R\rightarrow T$ so that $\phi$ has no ``edge folds''
(see below); then Skora \cite{Sk1,Sk2} showed that $\phi$ has no
``vertex folds'', giving that $\phi$ is a $\Gamma$-equivariant 
isometry, completing an
affirmative solution to the conjecture. 

\begin{theorem}[Morgan-Otal, Skora]
\label{theorem:skora}
Let $\Gamma=\pi_1(S)$, $S$ a closed surface of genus at least two.
Then any small, minimal $\Gamma$-action on an $\R$-tree is dual to the lift
of  a measured foliation on $S$.  
\end{theorem}

A complete exposition of Theorem \ref{theorem:skora} is given in \cite{Ot}.
\medskip

The purpose of the present paper is to prove Theorem
\ref{theorem:skora} from a different point of view, using harmonic maps.  
Harmonic maps were used by Gromov-Schoen \cite{GrSc} 
to show that certain groups do not act nontrivially on 
singular spaces such as trees.  Here we use harmonic maps to classify, in 
the special case of a surface group, all minimal, small 
actions on $\R$-trees, against a 
background where many such actions exist (namely $6g-7$ dimensions worth).

Our other interest in this proof is in the way it uses harmonic maps as a 
tool in combinatorial group theory.  For example, 
combinatorial topology arguments become greatly simplified (via the 
maximum principle) when looking 
at a harmonic representative.  Another example is the existence of a 
moduli space of harmonic maps (and harmonic map invariants) associated to 
a group action, allowing for an extra tool in proofs.

\subsection{Outline}
Here is a brief description of our approach to the proof. 

\medskip
\noindent
{\bf Step 1 (Find a harmonic map):} Given a small action 
of the surface group $\Gamma=\pi_1(S)$ on an $\R$-tree $T$, 
it is relatively straightforward to
find a $\Gamma$-equivariant harmonic map 
$f:\widetilde{S}\to T$. Here we have endowed $S$ with a complex
structure. 

\medskip
\noindent
{\bf Step 2 (Associated data): }
The harmonic map $f$ automatically has associated to it the following 
data:  
\begin{itemize}
\item a $\Gamma$-equivariant 
holomorphic quadratic (Hopf) 
differential $\widetilde\Phi$ on the Riemann surface
$\widetilde{S}$
\item a $\Gamma$-equivariant measured foliation $\widetilde{\cal F}$, 
the vertical foliation of $\widetilde{\Phi}$
\item the leaf space $R$ of $\widetilde{\cal F}$, with metric induced from 
the measure on $\widetilde{\cal F}$, making $R$ 
into an $\R$-tree
\end{itemize}

The map $f$ is projection along the
leaves of $\widetilde{\cal F}$, with the possibility 
of several vertical leaves being sent to the same point in $T$.  The 
$\Gamma$-action on $S$ induces an isometric $\Gamma$-action on $R$.

\medskip
\noindent
{\bf Step 3 (Morphism from a geometric action to the given action): }
Let $\pi:\widetilde{S}\to R$ be the natural projection 
sending each leaf of 
$\widetilde{F}$ to a point.  Here a leaf may have a countable 
number of k-pronged singularities.  
We then obtain a $\Gamma$-equivariant 
{\em morphism} $\phi:R\rightarrow T$ of $\R$-trees, where 
$\phi=f\circ\pi^{-1}$.   We must show that $\phi$ is an isometry, which 
is the same as saying that $\phi$ does not {\em fold} at any point.

\medskip
\noindent
{\bf Step 4a (No edge folds): }
If $\phi$ folded at an {\em edge point} of $R$, i.e. a point whose 
``tangent space'' has only two directions, then this would 
force $f$ to locally take the form $z\mapsto|\ree\ z|$
which is not harmonic.  Hence there are no edge folds, 
nor even {\em vertex folds} at trivalent vertices.  The {\em vertex 
points} of $R$ are precisely the images under $\pi$ of
leaves of  $\widetilde{\cal F}$ passing 
through a singularity of $\widetilde{\cal F}$.

\medskip
\noindent
{\bf Step 4b (No vertex folds): }
Ruling out folds at high order vertices $v\in R$ 
requires a global argument (see Example 3.2.1 of a 
local vertex fold).  The smallness hypothesis implies that, if two 
edges adjacent to $v$ are folded together, then neither edge can contain 
a point representing the lift of 
a closed leaf of $\cal F$.  
This basically allows us to reduce the proof to 
the model case (see \S 5.2.3) 
where some leaf of $\cal F$ is dense in $S$.  

We now exploit the fact that we have a choice of conformal 
structures for $S$.  Assuming 
$\phi$ folds at some vertex point, we can always choose a path of 
conformal structures $S_t$ on $S=S_0$ so that 
the Hubbard-Masur differential on $\widetilde{S}_t$ (the
holomorphic differential $\widetilde{\Psi}_t$ whose vertical 
foliation projects to $R$) has simple zeroes for $t\neq0$,
and the edges which are folded together are 
represented on $\widetilde{\cal F}_t$ by domains with
a common one-dimensional frontier. As the harmonic map would again
take the form $z\mapsto|\ree\ z|$ across this frontier,
we see that $\widetilde{\Psi}_t\neq\widetilde{\Phi}_t$, 
where $\widetilde{\Phi}_t$ is the Hopf differential for the harmonic map
$f:\widetilde{S}\to T$.  

Hence there is a family of distinct $\R$-trees $R_t$ and morphisms 
$\phi_t:R_t\to T$.  These
trees come from measured foliations on $S$ which themselves come from 
interval exchange maps. But any nontrivial continuous variation 
in an interval exchange
map forces a nontrivial variation in the tree $T$. As $T$ is fixed this is 
impossible, so there can be no vertex folds.

\subsection{Acknowledgements} 
We thank Mladen Bestvina for 
useful discussions, Howard Masur for all his help
(including the idea and most of the details of the
proof of Proposition 2.3), and the referee for
numerous comments and corrections 
which greatly improved the paper. 
Misha Kapovich \cite{Ka} also had the
idea to use harmonic maps
in the proof of Skora's theorem.

\section{Preliminaries}
\label{section:preliminaries}

\subsection{$\R$-trees}
 
An {\em $\R$-tree} is a metric space $T$ such that any two points in $T$ 
are joined by a unique arc and every arc is isometric to a segment in 
$\R$.  Let $[x,y]$ denote the unique (geodesic) segment from $x$ to
$y$ in $T$.

We say that $x\in T$ is an {\em edge point} (resp. {\em vertex point}) 
of $T$ if $T-\{x\}$ has precisely
two (resp. more than two) components.  An {\em edge} in $T$ is a nontrivial, 
embedded segment $[x,y]$ in $T$.
  
A {\em morphism} of $\R$-trees is a map $\phi:T\to T'$  such that 
every nondegenerate 
segment $[x,y]$ has a nondegenerate subsegment $[x,w]$ for which $\phi
\mid _{[x,w]}$ is an isometry.

The morphism $\phi:T\to T'$ {\em folds} at the point $x\in T$ if there
are nondegenerate segments $[x,y_1]$ and $[x,y_2]$, with 
$[x,y_1]\cap [x,y_2]=\{x\}$, such that $\phi$ maps each segment
$[x,y_i]$ isometrically onto a common segment in $T'$.  It is easy to
see that the morphism $\phi:T\to T'$ is an isometric embedding 
unless $\phi$ folds at some point $x\in T$.

By an {\em action} of $\Gamma$ on $T$ we mean an action by isometries.
The action is {\em minimal} if $\Gamma$ leaves no proper subtree of $T$
invariant.  For 
any $\Gamma$-action on $T$, there is
a $\Gamma$-invariant proper subtree which is minimal (see, e.g., \cite{CM}).
Also, if $A_\gamma$ is the isometry of $T$ corresponding
to $\gamma\in\pi_1\SR$ for which $\inf_{y\in T}d(A_{\gamma} y,y)>0$, then
$A_\gamma$ has an {\it axis\/} $l_\gamma$ in $T$, i.e., an isometrically
embedded line in $T$ which is invariant under $A_\gamma$ and which
has the property that $x\in l_\gamma$ iff $d(A_{\gamma} x,x)=\inf_{y\in T}
d(A_{\gamma} y,y)$. The proof is a straightforward consequence of the
non-positive curvature of $T$.
\medskip

\noindent
{\bf Assumption:} Henceforth we will assume, without loss of generality, 
that
all actions are minimal.
\bigskip

We will need the following fact about small actions.

\begin{lemma}
\label{lemma:benson} 
Let $\Gamma=\pi_1(S)$, $S$ a closed surface of genus at least two.  
If the action $\Gamma\x T\to T$ is small then $T$ must have a vertex
point.
\end{lemma}

Whenever speaking of vertex or edge points in a subtree 
of a tree $T$, we mean with respect to the space of directions 
in the subtree, not the ambient tree $T$.

\proof If $T$ has no vertex points then it is isometric to $\R$, so the
action of $\Gamma$ gives a representation $\psi:\Gamma\rightarrow
\mbox{Isom}(\R)$.  As $\psi(\Gamma)<\mbox{Isom}(\R)$ is 
virtually abelian and $\Gamma$ is
not solvable, it must be that the kernel of $\psi$ contains two
noncommuting elements.  But $S$ is a closed surface of genus at least 
two, so sufficiently high powers of any two noncommuting elements of
$\Gamma=\pi_1(S)$ generate a free group.  This free group lies in the 
kernel of
the action, in particular stabilizes any nondegenerate edge of $T$, a 
contradiction.
\endproof

\subsection{Holomorphic quadratic differentials}
\label{section:hqd}

A holomorphic quadratic differential $\Phi$ on the Riemann surface
$\SR$ is a tensor given locally by an expression $\Phi=\vp(z)dz^2$
where $z$ is a conformal coordinate on $\SR$ and $\vp(z)$ is
holomorphic. Such a quadratic differential $\Phi$ defines a measured
foliation in the following way. The zeros $\Phi^{-1}(0)$ of $\Phi$
are well-defined and discrete.  Away 
from these zeros, we can choose a canonical
conformal coordinate $\z(z)=\int^z\sqrt\Phi$ so that $\Phi=d\z^2$.
The local measured foliations ($\{\ree\z={\mbox{const}}\}$,
$|d\ree\z|$) then piece together to form a measured foliation known
as the vertical measured foliation of $\Phi$.

\subsection{Actions dual to a measured foliation}
\label{section:dualaction}

Let $(\SF,\mu)$ denote the vertical
measured foliation of $\Phi$.  Lift it to a $\pi_1\SR$-equivariant
measured foliation $(\wt\SF,\tl\mu)$ on $\wt\SR$. The leaf space $R$
of $\wt\SF$ is a Hausdorff topological space.  Let 
$\pi:\wt\SR\to R$ denote the projection.  The leaf space 
$R$ of the measured foliation
$(\widetilde{\cal F},\mu)$ inherits a
metric space structure from the measure $\mu$: a geodesic segment $[x,y]$ 
in $R$ is given by any path $\gamma$ in $\hyp^2$ from a point in 
the leaf corresponding to $x$ to a point in the leaf corresponding to $y$, 
such that $\gamma$ is transverse to the leaves of the foliation $\widetilde{\cal
F}$. The distance $d_R(x,y)$ is simply $\mu(\gamma)$, and the metric space
$(R,d)$ is an $\R$-tree (see \cite{MS2}).  
This tree is often not locally compact.  For instance, when 
the leaves of the foliation on the surface $S$ are dense, we can find 
sequences of arcs $C_n$ transverse to
the foliation with endpoints on singularities
of $\widetilde{\cal F}$ whose transverse measure $\mu(C_n)$ goes to zero,
forcing the distance between the corresponding images of the (lifts
of) vertices to also go to zero.  

The action of $\Gamma$ on $\hyp^2$ preserves $\mu$, and so induces an
isometric action of $\Gamma$ on $R$.  The stabilizers of this action are
virtually cyclic, in particular are small.   

The action of $\pi_1\SR$ on $\hyp^2$ preserves $\mu$, and so induces an
isometric action of $\pi_1\SR$ on $R$.  
The map $\pi:\wt\SR\to R$ is equivariant with respect
to this action.

\subsection{The Hubbard-Masur Theorem}

Holomorphic quadratic differentials on a Riemann surface $\SR$ are related 
to
classes of measured foliations via the Hubbard-Masur Theorem. To set the
notation, fix a Riemann surface $\SR$ and define a map $HM:\qd(S)\to 
\MF(S)$ from
the space $\qd(S)$ of holomorphic quadratic differentials on $\SR$ to the
space $\MF(S)$ of equivalence classes of measured foliations on $\SR$ that
associates to  $\Phi\in\qd(S)$ the class of its vertical measured
foliation. A fundamental result is

\begin{theorem} [Hubbard-Masur \cite{HM}] \label{theorem:Hubbard-Masur} 
\label{theorem:hm}
HM is a surjective homeomorphism.
\end{theorem}

\remark A proof of Theorem \ref{theorem:hm} 
in the spirit of the current work can be found in
\cite{W2}.
\bigskip

An alternative phrasing will be convenient for us. Let $Q(S)\to\Teich(S)$ 
denote the bundle of holomorphic quadratic 
differentials over $\Teich(S)$: here the
fiber over $[S]\in\Teich(S)$ is the space $\qd(S)$ of quadratic 
differentials holomorphic with respect to a complex structure 
$S$ in $[S]$. Let $(\SF,\mu)$ denote
a given measured foliation. Then the Hubbard-Masur Theorem shows that 
there is a
well-defined section $\Psi_\mu:\Teich(S)\to Q(S)$ which associates to
$[S]\in\Teich(S)$ the holomorphic quadratic differential 
$\Psi_\mu (S)\in \qd(S)$ whose
vertical measured foliation is measure equivalent to $(\SF,\mu)$.

\subsection{Moving around in the Hubbard-Masur section}

In this subsection we give a basic property of the section $\Psi_\mu$.

Let $S$ be a Riemann surface and let $q$ be a holomorphic
quadratic differential with vertical measured foliation 
$(\SF,\mu)$. Let $p_0$ be a singularity of $q$ and let $L$
be the maximal compact graph of singular vertical arcs
through $p_0$ which connect $p_0$ to all the other singularities
on the leaf through $p_0$.
Consider a neighborhood $\cal N$ of L in $S$.  We refer to the
components $\{s_i\}$ of ${\cal N} - L$ as {\it sectors}, and say that 
two sectors {\em meet} along a (nondegenerate) 
arc if their closures intersect along 
that arc.  We observe that there is a natural correspondence of sectors
near a maximal singular arc $L$ as above under Whitehead moves
and isotopies of the foliation.

\begin{proposition}[Sectors can be made adjacent]
\label{proposition:moving}
Let $S$ be a Riemann surface, let $q$ be a holomorphic
quadratic differential with vertical measured foliation 
$(\SF,\mu)$, and let $L$, $p_0$ and  $\{s_i\}$ be as above.
Choose any pair of sectors 
$s_{i_1}$ and $s_{i_2}$ from the list of sectors  $\{s_i\}$.       
Then there is a Riemann surface $\SR^*$ and a
holomorphic quadratic differential $q^*$ on $\SR^*$ so that the vertical
foliation of $q^*$ is measure equivalent to $(\SF,\mu)$, and under the
equivalence the sectors $s_{i_1}^*$ and $s_{i_2}^*$ (corresponding
to $s_{i_1}$ and $s_{i_2}$, respectively) meet along an arc.
\end{proposition}

Note that both $q$ and $q^*$ are in the image of the Hubbard-Masur 
section corresponding to $(\SF, \mu)$.  A self-contained
proof of Proposition 
\ref{proposition:moving} is given in the appendix.

\section{Harmonic maps to trees}

\subsection{Definition of harmonic map}
\label{section:harmonicmaps}

Given a Lipschitz continuous map $w: \SR\to(T,h)$ from a 
Riemann surface $\SR$ to a locally finite metric tree $T$, we
define the energy form to be the tensor
$$edz\ox d\bar z=(\|w_*\partial_z\|_h^2+\|w_*\partial_{\bar z}\|_h^2)dz\ox
d\bar z$$
Since the map $w$ is Lipschitz, it is differentiable almost
everywhere and bounded almost everywhere on closed balls; thus the
form $edz\ox d\bar z$ is defined almost everywhere with 
$edz\wedge d\bar z$ integrable over compacta.  Note that even when $T$ is 
not locally finite, the image of any closed ball in $S$ is compact hence 
lies in a locally finite subtree of $T$, so this analysis applies.

Alternatively, for any conformal metric $g$ on $S$ with 
area form $dA_g$, the energy form may be expressed as follows.
Choose an orthonormal frame $\{v_1,v_2\}$ at a point $z\in S$,
and consider the pushforward vectors $\{w_*v_1, w_*v_2\}$. The the
energy form is the 2-form 
$\frac12(\|w_*v_1|_h^2+\|w_*v_2\|_h^2)dA_g$,
or alternatively $\frac12tr_g(w^*h)dA_g$.
The energy of the map $w$ is $E=\int edz\land d\bar z$. The map $w$ is 
a {\em harmonic map} if it is a minimum for this functional in its 
homotopy class of maps. We define the Hopf differential $\Phi$ for a map
$w: \SR\to T$ by
$$\Phi=\phi dz^2=4\left<w_*\partial_z, w_*\partial_z\right>_hdz^2$$
 
Note that $\|\Phi\|=\|\Phi\|_{L^1}<2E$.

\subsection{Examples}

In this subsection, we list some motivating examples of harmonic maps from
Riemann surfaces to $\R$-trees. Each example will illustrate a
principle we will later use.

\begin{enumerate}
\item The map $f(z)=\ree\{z^2\}$ as the most basic vertex fold.

The map $f(z)=\ree\{z^2\}:\C\to\R$ can be viewed as a harmonic map from the
Riemann surface $\C$ to the $\R$-tree $\R$. Observe that the preimage of 
the
origin $O\in\R$ is the pair of intersecting lines $\{x=\pm y\}$ which 
divides
$\C$ into four sectors. The other level lines of a nonzero $r\in\R$ 
consist of
hyperbolas $\{x^2-y^2=r\}$.  The leaf space of the {\em connected\/} 
components
of these level curves is the pair of coordinates axes. We conclude that the
harmonic function $f(z)$ factors as a projection to the $\R$-tree of the
coordinates axes followed by a vertex fold of each half-axis to its 
negative,
which results in the image $\R$-tree $\R$.

\item Here is an example from \cite{W3}:
begin with the holomorphic differential $z^kdz^2$ on $\C$, whose
vertical measured foliation is the set of curves $\{Re z^{k+2} = c\}$.
When we project along this foliation, we obtain a harmonic map to a 
tree with $k+2$ prongs out of a single vertex. 

\item Actions Dual to a Measured Foliation $({\cal F},\mu)$, as given
  in \S\ref{section:dualaction}.

Here the harmonic map is simply the projection along the vertical 
foliation of
the properly normalized Hubbard-Masur differential for $(\SF,\mu)$.  
This characterization is independent of the particular 
Riemann surface chosen. We therefore observe the following.

\begin{lemma}[When $\Phi$ and $\Psi$ agree]
\label{lemma:3.1} 
If the action of $\Gamma$ is dual to a
measured foliation $(\SF,\mu)$, then there is a well-defined
Hopf differential 
section $\Phi:\Teich(S)\to Q(S)$ for $\pi$, and this section
$\Phi$ agrees with the 
Hubbard-Masur differential section $\Psi:\Teich(S)\to Q(S)$
for $\SF$.
\end{lemma}

\proof The lemma is effectively the content of [W2], which we now
summarize; for complete details, see [W2].  (Later on,
in \S4, we shall give an independent proof of the existence 
of a harmonic map dual to a measured foliation.)
As in \S2.3, a measured foliation $(\SF,\mu)$ on $S$ 
lifts to an equivariant measured foliation $(\wt\SF,\tl\mu)$
on $\tl S$; we can project along the leaves to obtain an
$\R$-tree $(R, d)$, with this construction also yielding an
equivariant map  $\pi_0: \tl S \to (R,d)$. 

For each complex
structure $\sigma$ on $S$, we can minimize the energy in the 
equivariant homotopy class of $\pi_0$ obtaining [W2; proof
of Prop. 3.1] a map $\pi: (S,\sigma) \to (R,d)$ whose Hopf
differential $\Phi_{R}(\sigma)$ has vertical foliation 
measure equivalent to $(\SF,\mu)$. (This argument is a straightforward
application of Ascoli-Arzela, with a crucial use of the axes of
group elements in $R$ to control (see [W2; Lemma 3.4]) the images
of some points by elements of the minimizing sequence of maps.)
This characterizes the differential uniquely \cite{HM}; for a harmonic
maps argument for this uniqueness, see [W2; \S4]. Here the point is 
that both maps can be given as projections along 
minimal stretch foliations of Hopf differentials and 
the distance between the image points of the two maps 
can be equivariantly defined, and is a subharmonic function.
[As the pullback of a smooth convex function off of the zeroes of the
Hopf differentials, this pullback of the distance function is 
smooth and subharmonic (i.e. submean for balls of fixed radii)
away from a discrete set of singularities and continuous
across them; hence it is subharmonic everywhere. (Compare 
Proposition~3.2)] The maximum
principle then applies, showing that the distance must be
constant. Some geometry of the
tree, in particular the fact that it has branches, forces that
constant to vanish. Thus $\Phi_{R}(\sigma) = \Psi_{\mu}(\sigma)$. 
\endproof

An important part 
of our proof of Theorem \ref{theorem:skora} will be a converse 
(Lemma \ref{lemma:key}) to 
Lemma \ref{lemma:3.1}.

\item Actions Dual to the Measured Foliation of the Hopf 
Differential for an Arbitrary Harmonic Map from a Surface.

Let $f:S \to X$ be a harmonic map from the Riemann surface $S$
to a metric space (possibly Riemannian, possibly singular).
Let $\Phi$ denote the associated Hopf differential; we will see
in \S\ref{section:localstructure}, that this Hopf differential
is a holomorphic quadratic differential on $S$.  Lift this 
differential $\Phi$ to an (equivariant) differential 
$\widetilde{\Phi}$ on the universal cover $\widetilde{S}$,
and consider its vertical (corresponding to the minimal 
stretch directions of $\tilde{f}$) measured foliation
$(\widetilde{\cal F},\mu)$ and associated projection
$\pi:\widetilde{S} \to R$ to the leaf space 
(see \S\ref{section:dualaction}).

Part of the content of the previous example is that the equivariant
projection map $\pi:\widetilde{S} \to R$ from $\widetilde{S}$ to
$R$ is harmonic.

Our proof of Skora's theorem involves a study of the relationship
between the harmonic map we will construct from $\widetilde{S}$ to $T$ and the
associated harmonic map $\pi:\widetilde{S} \to R$ from $\widetilde{S}$
to the leaf space $R$.

\remark In [W4], we study the product
harmonic map $(\tilde{f},\pi):\widetilde{S} \to X\times \R$, and find that 
it is also
conformal, after a slight homothety of $\R$. We also find that when
$X$ is smooth and two-dimensional, then this map is a stable
minimal map.

\item A Harmonic One-Form with Integral Periods.

Project a square torus $T^2$ along its natural vertical foliation to the 
circle
$S^1$. This map is clearly harmonic. Now, there is a genus two surface 
$\SR$
which is a branched cover over $T^2$, and the $1$-form $dz$ lifts to
yield a holomorphic one form on $S$. One can still project along 
leaves of this
one form to a figure~8 which is a branched cover of the original
$S^1$.  Hence by composing with the map to $S^1$, we see 
that there is an associated harmonic
map $f:\SR\to S^1$ and, via the one-form lifted from $dz$, an associated
holomorphic one-form with integral periods on $\SR$. 

As we vary $\SR$ in
Teich~$(\SR)$ (say in a family $\SR_t$), the holomorphic one forms with the
same $A$-periods (in the usual notation) varies continuously through 
one-forms,
say, $\om_t$. It is interesting to consider the topology of the foliations
$\SF_t$ that integrate $\ker\ree\om_t$. 

The original surface $\SR$ could be described as being constructed from a 
pair
of cylinders $C_1$ and $C_2$ bounded by circles $\{S_{11},S_{12}\}$ and
$\{S_{21},S_{22}\}$, respectively. Each $S_{ij}$ is composed of two 
semicircles
$s^{\mbox{t}}_{ij}$ and $s^{\mbox{b}}_{ij}$. Now, the upshot of the
notation is that $\SR$ is defined by identifying $s^{\mbox{b}}_{11}$ to
$s^{\mbox{t}}_{12}$, $s^{\mbox{b}}_{12}$ to $s^{\mbox{t}}_{21}$,
$s^{\mbox{b}}_{21}$ to $s^{\mbox{t}}_{22}$ and $s^{\mbox{b}}_{22}$
to $s^{\mbox{t}}_{11}$ in the natural way, and the foliations parallel to 
the
core curves of the cylinders become $\SF_0$.

A natural motion in \tec space is to slightly rotate one of these cylinders
against the other. This has the effect in our case of preserving the 
topology
of $\SF_0$, up to a Whitehead move which splits the singularities from 
being
locally a pair of coordinate axes as in Example~1 to a ``double $Y$''
configuration. Of course, as $\SF_t$ is the foliation of a harmonic 
one-form,
we see that this new synthetically constructed foliation inherited from the
cylinders, which is actually the foliation of the Hubbard-Masur 
differential
for $(\SF_0,\mu_0)$ on $\SR_t$, is not $\SF_t$. We conclude in this case 
that
for the natural $\pi_1S$- action on $\R$ defined via the one form $\om$, it
is not the case that the Hopf differential section 
$\Phi:\Teich(S)\to Q(S)$ agrees with
the Hubbard-Masur differential section $\Psi:\Teich(S)\to Q(S)$.
\end{enumerate}

\subsection{Local structure}
\label{section:localstructure}

R. Schoen has emphasized (see \cite{Sc}) 
that a map for which the energy functional
is stationary under reparametrizations of the domain has a Hopf
differential which is holomorphic: one uses suitable domain
reparametrizations to show that the Hopf differential satisfies the
Cauchy-Riemann equations weakly, and then Weyl's Lemma forces the
Hopf differential to be (strongly) holomorphic. We observe that in
this argument, the range manifold may be singular.

The vertical and horizontal foliations of the Hopf differential for
$w: \SR\to T$ integrate the directions of minimal and maximal
stretch of the gradient map $dw$, for smooth energy minimizing maps
$w: \SR\to T$. As the image is one-dimensional, the harmonic map $w$ is
a projection along the minimal stretch direction.  Further,
if one normalizes the conformal coordinates in a domain that
avoids the zeroes $\Phi^{-1}(0)$ of the Hopf differential $\Phi$
so that $\Phi = dz^2$ in that neighborhood, then one sees
from the geometric definition of $\Phi$ above that the 
energy-minimizing map takes maximal stretch segments of measure
$\epsilon$ to segments in $T$ of length $\epsilon$.

\subsection{Effect on convex functions}

A function defined on a an $\R$-tree is {\em convex} if its restriction 
to every geodesic is convex in the classical sense.  Recall that a 
function is {\em subharmonic} if it is submean, that is it's value at any 
point is less than or equal to its average in a small ball around that 
point.  Harmonic maps between Riemannian manifolds pull back convex 
functions to subharmonic functions (see, e.g.\, \cite{I}).  This 
important property extends to the case of $\R$-tree targets.

\begin{proposition}
A harmonic map from a surface $S$ to an $\R$-tree pulls back
(germs of) convex functions to (germs of) subharmonic functions.
\end{proposition}

\proof
We first argue that the map $\pi: \widetilde{S} \to R$ to the leaf space 
$R$ pulls back germs of 
convex functions on $R$ to germs of subharmonic 
functions on $\widetilde{S}$.
Locally, the level set of the vertex, say  $V \in R$, 
near a zero or pole of the Hopf differential
$\widetilde{\Phi}$
divides the neighborhood of the singularity
into 'sectors', with the natural coordinate $\zeta$ mapping each 
sector conformally onto a neighborhood of zero 
in the upper half-plane  
(see [St;\S7.1]).  Under this mapping of 
a sector, the foliation
of preimages of points in the tree $R$ (in a sector)
is taken to the horizontal foliation of the 
half-plane given by curves of the form $\{y=const\}$.

While it is not essential for the proof at this point,
let us now
consider a convex function $F$ 
defined on the tree $R$ near the point 
$V \in R$. This function pulls back to a function on a collection
of sectors, which is constant on the horizontal levels
 $\{y=const\}$, and convex in $y$. Since any sector 
can be taken to any other sector by an appropriate rotation,
it is straightforward to see that this pullback is submean.
(A more detailed argument is also given below, in the case
of the tree $T$.)

With this in mind, let us return to the original case of 
the map $f:\widetilde{S} \to T$. In the neighborhood of the
singularity $p$ of the Hopf differential $\widetilde(\Phi)$, 
we can regard our map as first projecting to a neighborhood
of a vertex $V$ in R (this neighborhood of $V$ is metrically
a k-pronged star out of $V$, with one prong for each sector,
by construction), followed by a map of $R$ to $T$, in which 
several prongs of $R$ map to a single prong of
$T$, this prong of T emanating out of the
image $v \in T$ of the vertex $V \in R$). Here
we must have each prong taken injectively to a prong,
because of the form of the map $f:\zeta \mapsto Re \zeta = \xi$
away from singularities of $ \widetilde(\Phi)$.

In order to see why the map $f:\widetilde{S} \to T$ 
pulls back germs of 
convex functions on $T$ to germs of subharmonic 
functions on $\widetilde{S}$, we make one crucial 
observation: we note that neighboring sectors 
on $\widetilde{S}$ must
be taken to different prongs out of $f(p) \in T$;
this is because a small arc transverse to the common 
boundary leaf of the pair of sectors is projected
by $f$ injectively into $T$, once again because 
in a neighborhood of such an arc, there are no
singular points, and so the map $f$ is of the form
$\zeta \mapsto Re \zeta = \xi$.
This implies that the pre-image under $f$ of any 
given prong in $T$ consists of at most half of 
the sectors abutting $p$.

Consider then  a convex function $F$ on the 
tree $T$ near a point 
$p \in T$. This function pulls back to a function 
on a collection
of sectors, which is constant on the horizontal levels
$\{y=const\}$, and also convex in $y$.  Suppose we have
that $F(v) =0$, so we need the mean value of $f^*F$
to be non-negative on a disk $D$ around $p$.  Of course,
since $F$ is convex on $f(D) \subset T$, we know that
$F$ can be negative on at most one prong of $f(D)$, and
must be non-negative on the other prongs. Moreover,
since $F$ is convex, if we average $f^*F$ over a pair
of sectors, one in which $f^*F$ is nonpositive, and
one in which $f^*F$ is non-negative, we see that the
sum of the averages must be non-negative. (To see this,
conformally map each sector to a half-plane
(say $\{ y\ge0\}$), and then
glue the halfplanes together so that $f^*F$ is convex
in the coordinate $y$ across the foliation.) Then we simply
apply the observation of the previous paragraph to conclude
that since $f^*F$ is non-negative on at least half of the sectors,
the average of $f^*F$ on the union of sectors (i.e th 
disk $D$) must be non-negative, as required. 
\endproof

\section{Constructing a morphism from a geometric action to the given
  action} 
\label{section:construction}

Let $\Gamma=\pi_1(S)$, $S$ a closed surface of genus at least 2, and let
$\Gamma\times T\to T$ be an action (not necessarily small) on
an $\R$-tree $T$.  In this section we 
construct an action of $\Gamma$ on an $\R$-tree $R$
which is dual to a measured foliation, and a $\Gamma$-equivariant 
morphism $\phi:R\to T$.  

We will think of $S$ as having a fixed hyperbolic
structure, and so the universal cover $\widetilde{S}$ is the
hyperbolic plane $\hyp^2$.  Since $T$ is contractible, there is a
$\Gamma$-equivariant Lipschitz 
continuous map $f_0:\hyp^2\rightarrow T$.  To be
concrete, one may lift a triangulation of $S$, 
define $f_0$ by equivariance on the 0-skeleton of this triangulation,
then extend (by contractibility of $T$) equivariantly to the
1-skeleton and 2-skeleton.

\subsection{Finding the foliation using a harmonic map}

Our first goal is to find a harmonic $f$ map in the 
equivariant homotopy class of the $\Gamma$-equivariant 
continuous map $f_0:\hyp^2\rightarrow T$ constructed at the beginning
of \S\ref{section:construction}.  The harmonic
map $f$ will have the property that there is a measured foliation 
$({\cal F},\mu)$ on $S$ so that every leaf of $\widetilde{\cal F}$
gets mapped to a point under $f$.  While it is possible to use the 
general
theories of Korevaar-Schoen \cite{KS} and Jost \cite{J1,J2} 
on harmonic maps to nonpositively curved metric spaces, we will
construct the harmonic map from elementary methods here.

To carry this out, we choose balls $B_1,\ldots , B_n$ on $S$ so that:
\begin{itemize}
\item the balls are topologically trivial 
\item the restriction
$f_0\mid _{\hat{B}_j}$ of $f_0$ to a lift $\hat{B}_j$ of $B_j$ is not a
constant map for $j=1,\dots,n$, and 
\item the set $\{B_1,\ldots B_n\}$ of
balls is an open cover of $S$
\end{itemize}

Thus we have that each lift of
$B_j$ is disjoint from every other lift of $B_j$, and the union of 
all the lifts of all the balls $\{B_1,\ldots B_n\}$
covers $\widetilde{S}$.

Now for each lift $\widehat{B}_1$ of $B_1$ the image $f_0(\widehat{B}_1)$ 
is a
finite subtree of $T$.  This follows from the fact that, for a
basepoint $b_1\in \widehat{B}_1$, the image $f_0(\widehat{B}_1)$ lies 
in a compact subset $K$ of the space of directions at $f_0(b_1)$, 
and as this space of directions $K$ is discrete (from the definition of 
$\R$-tree), it is also finite. 

It is straightforward that there exists a
unique harmonic map $\widehat{f_0}:\widehat{B}_1\to T$ so that
$\widehat{f_0}\mid _{\partial\hat{B}_1}=f_0\mid_{\partial\hat{B}_1}$ 
(see the Appendix of \cite{W1} for existence.
To see uniqueness, apply the method of Cor. 3.2 of [W3]
(see also \S4 of [W1]): the distance between any pair of solutions
would be subharmonic on $\widehat{B}_1$ and
vanishing on $\partial\hat{B}_1$ -- thus any pair of solutions
coincide.). Moreover, if
$h(\widehat{B}_1)$ is any other lift of $B_1$, the uniqueness of the
harmonic map then forces 
$\widehat{f_0}\mid_{h(\hat{B}_1)}=h\circ\widehat{f_0}\mid_{\hat{B}_1}$.  
Let $\phi_1$ denote the map from the complete lift of $B_1$ to $T$. 
Then $\phi_1$, being nonconstant, also has the following properties:
\begin{itemize}
\item $\phi_1$ is projection along the vertical 
measured foliation of its Hopf
differential, and 
\item $\phi_1$ is  $C^\infty$ on the interior of its
domain (off of the zeroes of the Hopf differential of $\phi_1$)
\end{itemize}

Set 
$$
f_1(z)=
\left\{
\begin{array}{ll}
\phi_1(z) & z\in 
\mbox{\ lift of $B_1$} \\
f_0 & \mbox{otherwise}
\end{array}
\right.
$$

Then $f_1$ is equivariantly homotopic to $f_0$,
and is a $C^\infty$ equivariant projection (as above) along a measured
foliation on the domain of $\phi_1$.

We repeat this procedure for lifts of the ball $B_2$, using $f_1$
as the original map instead of $f_0$. We then obtain a map $f_2$. The
situation is most interesting when $B_1\cap B_2\neq\emptyset$, as
then the boundary values for $\phi_2$ are defined by values of
$\phi_1$, which may not agree with those of the original $f_0$. 

The
main observation we need to make is the following: along most of a
small neighborhood of $\partial \widehat{B}_2\subset \widehat{B}_1$ 
we have that 
$\phi_1 \mid _{\hat{B}_1\setminus \hat{B}_2}$ and $\phi_2 \mid
_{\hat{B}_2}$ extend to 
be a well-defined Lipschitz projection along a well-defined
Lipschitz measured foliation. To see this note that
$\phi_1 \mid _{\partial \hat{B}_1}$ is $C^{1,\alpha}$ and the measure of 
the
vertical foliation of the Hopf differential of $\phi_1$ is defined
by distance between image points in $T$ (see \S2.5). As this
also holds for $\phi_2 \mid _{B_2}$, and
$\partial \overline{\widehat{B}_1}\subset \overline{\widehat{B}_2}$ 
is compact, the claim follows, except
at (a discrete set of) places where the boundary values
$f_1 \mid_{\partial \overline{B_2}}$ double back and result in small arcs 
in
both $\overline{B_1}$ and $\overline{B_2}$ 
which close up in $\widehat{B}_1\cup \widehat{B}_2$.

We follow the same procedures iteratively for lifts of
$B_3,\dots,B_n$ obtaining an equivariant map $f_n:\widetilde{S}\to T$ which
is a projection along a Lipschitz measured foliation except for a
discrete set of places where the leaves are closed and
homotopically trivial.

In these places, we do an equivariant surgery to the map. For any
region consisting of a union of concentric closed leaves, consider
the closure of the largest such region. 
We then collapse the region to a segment which
maps to the point defined by the boundary leaves.  
Call the new (collapsed) map $F:\widetilde{S}\to T$. It is evidently an
equivariant map along a measured foliation with singularities that
are $k$-pronged.

In \cite{W2}(Prop. 3.1),
an elementary proof shows that the piecewise harmonic map
$F:\hyp^2\rightarrow T$ as above is equivariantly homotopic to a 
harmonic map $f:\hyp^2\rightarrow T$. (This proof only requires that there 
are two elements of $\Gamma$ whose axes in $T$ 
have unbounded intersection. This property is
much weaker than requiring that the action be small, but, for our purposes,
follows from Lemma \ref{lemma:benson} above.)   
Moreover, attached to $f$ is a holomorphic
quadratic differential $\wt\Phi_0$, the Hopf differential of
$f$, with the following  properties (see [W2; \S2.2]): 
\begin{itemize}
\item The vertical measured foliation of $\wt\Phi_0$ is 
equivalent to $(\wt\SF,\tl\mu)$. 
\item The leaf space of the vertical foliation
of $\wt\Phi_0$ is $R$, and the vertical measure pushes down (say via
$\pi:\hyp^2\to R$) to the metric on $R$. This map is harmonic.
\item On neighborhoods $B \subset \tl S$ which are disjoint from 
$\widetilde\Phi^{-1}(0)$, the map $f|_B$ agrees with $\pi|_B$ up to
an isometry, while $\pi|_B$ is the projection $z \mapsto Re z$
in a natural coordinate system.
\end{itemize}

This last property is quite important for the sequel, so we recall
some the details from, for instance, [W1; p. 273] and [W2; p. 117].
By \S\ref{section:hqd}, there is a canonical 
coordinate $\zeta = \xi + i\eta$
so that $\Phi_0 = d\zeta^2$ on $B$. In its guise as a Hopf
differential, of course, the definitions from
\S\ref{section:harmonicmaps} provides that 
$\Phi_0 = \parallel f_*\partial\xi \parallel^2 -
\parallel f_*\partial\eta \parallel^2 + 2i<f_*\partial\xi,
f_*\partial\eta>$. Combining these two descriptions of $\Phi_0$
and using that B is one-dimensional, we find that $f|_B$ is 
isometric to the map $\zeta \mapsto Re \zeta = \xi$.

\subsection{Definition of the morphism $\phi: R\to T$}

Define an associated harmonic projection
$\pi:\widetilde{S}\to R$ via the construction in Example 3.2.4.
Define also a map $\phi: R \to T$ by $\phi=f \circ \pi^{-1}$.
We claim that $\phi$ is a morphism. To see this, let $I$ denote a
nondegenerate segment on the tree $R$; we must find a non-degenerate
subsegment 
$J\subset I$ for which $\phi\bigm|_J$ is an isometry. Well, as $R$ is 
defined via  projection $\pi:\widetilde{S}\to R$, we can 
find an arc $\gamma\subset\widetilde{S}$ with $\pi(\gamma)=I$.  
Here $\gamma$ is quasi-transverse to $\widetilde{\SF}$ 
(in the sense of \cite{HM},p. 231) and
$\mu(\gamma)=\ell_R(I)$. On any subarc $\gamma'$ of $\gamma$ which avoids  
the zeros of $\wt\Phi_0$, we may write  
(as we did at the end of the previous subsection)
$\wt\Phi_0=dz^2$ for a choice of conformal
coordinate in a neighborhood of $\gamma'$, and 
(again as in the previous subsection) $f$ is an isometric 
submersion.
Then for $J=\pi(\gamma')\subset\pi(\gamma)=I$, we have that
$\phi\bigm|_J=f\bigm|_{\gamma'}$ which is an isometry by construction. 

Finally, $\phi$ is
surjective by the minimality hypothesis,  and $\phi$ is equivariant since
$f$ is equivariant.

\section{Proving that $\phi$ doesn't fold}
\label{section:nofolds}

\subsection{No edge folds}

It is a direct consequence of harmonicity of $\phi$ that $\phi$ does
not fold at edge points.  This is actually implicit in the proof above 
that $\phi$ is a morphism, but we give a slightly different proof in the 
next proposition, to which we will refer back several times in the sequel.

\begin{proposition}[no edge-point folds]
\label{proposition:edgefolds}
The morphism $\phi:R\to T$ does not fold at an edge point $x\in R$.
\end{proposition}

\proof
The pre-image of an edge point is a nonsingular leaf of the  
foliation $\widetilde{\cal F}$.  
Any point $z_0$ on this leaf has a neighborhood $\SN$ foliated
by non-singular arcs of leaves, and admits a conformal coordinate $z=x+iy$
with the foliation parallel to $\ker(\ree dz)$. If $\phi:R\to T$ were to 
fold at
an edge point $\pi(z_0)$, then the harmonic map on the neighborhood $\SN$ 
would necessarily have the form
$z\mapsto|\ree z|$, which is, of course, not harmonic. 

Alternatively, using the same notation for the 
morphism $\phi$ folding at an edge point $\pi(z_0)$, letting $p_0$
denote the point $p_0=\phi\circ \pi(z_0)$, we may 
apply the maximum principle to the function $h=f^*(-d_T(p_0,\cd))$ on a
neighborhood of $z_0$. Here $-d_T(p_0,\cd)$ is convex on $f(\SN)$, while
$f^*(-d_T(p_0,\cd))$ is not subharmonic on $\SN$, contradicting
Proposition 3.2.
\endproof

Note that, at this point, we have shown that for any small 
$\Gamma$-action on an $\R$-tree $T$, there is an action on a
tree $R$, dual to a measured foliation, and a $\Gamma$-equivariant
morphism $\phi:R\rightarrow T$ which folds only at vertex points.

\subsection{No vertex folds}
In this section we will show that, when
the action of $\Gamma$ on $T$ is small, the morphism 
$\phi$ is an isometry. A crucial
feature of our argument will be a lemma that says that for actions 
$\Gamma\x T\to T$
which are not small, the choice of tree $R$ is not uniquely determined.

\begin{proposition}[no vertex folds]
\label{proposition:vertexfolds}
If the action $\Gamma\x T\to T$ is small, 
then the morphism $\phi$ does not fold at a vertex point $v\in R$.
\end{proposition}

The rest of this section is devoted to proving Proposition 
\ref{proposition:vertexfolds}. 

\subsubsection{Vertex fold gives bad family}
\label{subsection:family}

We begin with the following generalization
of Lemma~\ref{lemma:3.1}.

\begin{lemma} 
\label{lemma:key}
With notation as above, the following conditions are equivalent:
\begin{enumerate}
\item The action of $\Gamma$ on $T$ is dual to the measured foliation 
$\cal F$.
\item The morphism $\phi:R\to T$ is an isometry.
\item The Hubbard-Masur section $\psi_{\cal F}:\Teich(S)\to Q(S)$ for 
$\cal F$ is the same as the Hopf differential section 
$\Phi:\Teich(S)\to Q(S)$ for $T$.
\end{enumerate}
\end{lemma}

\proof As $R$ is the dual tree of $\widetilde{\cal F}$, 
it is clear that (2) implies (1).  Lemma~\ref{lemma:3.1} states 
that (1) implies (3).  

Now we prove that (3) implies (2).  
If $\phi$ is not an isometry, then $\phi$ must fold at some vertex point 
$v\in R$, by  Proposition \ref{proposition:edgefolds}.  Say 
$\phi(e_1)=\phi(e_2)$ for edges $e_1,e_2$.

We may assume that $R$ has vertices of valence at least 4: 
otherwise a vertex fold at a vertex $v\in R$ would be the fold of a
$3$-pronged
star to an interval or half-interval. Thus the map $f$ would restrict, in a
neighborhood of the pre-image of $v\in R$, to a harmonic function
on a disk where the corresponding Hopf differential has a 3-pronged zero.  
This is impossible, as harmonic functions are
locally $\ree(cz^k)+O(z^{k+1})$, for $k$ an integer.

[Alternatively, the preimage of $\pi^{-1}(v)$ 
is a tree with discrete trivalent
singularities. Near the singularities, this 
tree locally disconnects $\wt S$ into
three sectors, with the harmonic map 
$f:\wt S\to T$ folding the image of one sector
onto the image of an adjacent sector. 
Yet as the sectors meet along an edge, the
proof of Proposition \ref{proposition:edgefolds} 
applies to yield a contradiction.]

Now consider the Hubbard-Masur section
$\psi_{\cal F}:\Teich(S)\to Q(S)$ 
for the foliation $(\SF,\mu)$.  
We are assuming that $\psi_{\cal F}(S)$ has zeroes of order at least two.  
Let $s_1,s_2$ be the sectors of $\cal F$ corresponding to the edges 
$e_1,e_2$.  By Proposition \ref{proposition:moving} there is another 
quadratic differental $q'=\psi_{\cal F}(S')$ so that the 
sectors of the vertical foliation of $q'$ corresponding to $s_1,s_2$ have 
closures which meet along an edge.  Since by assumption $\psi_{\cal F}$ 
is the same as the Hopf differential section $\Phi$, this is a 
contradiction: it violates the maximum principle for the map $f'$
(defined as projection along the foliation of 
$\psi_{\cal F}(S') = \Phi(S'$)), again as the
map would locally have the form $z\mapsto|\ree\ z|$. 
\endproof

\bigskip
\noindent
{\bf Proof of Proposition \ref{proposition:vertexfolds}: }
We now suppose, in expectation of reaching a contradiction, that the 
given action is not dual to a measured foliation, i.e. that $\phi$ 
is not an isometry.  The equivalence of (1) and
(3) in Lemma \ref{lemma:key} then implies that there is a family
$\{S_t\}, t\in \R$ of distinct Riemann surfaces
for which $\psi_{\cal F}(S_t) \ne \Phi(S_t)$ 
for $t>0$ and $\psi_{\cal F}(S_0)= \Phi(S_0)$
(here $\cal F$ is defined by setting $\psi_{\cal F}(S_0)= \Phi(S_0)$
for some base point $S_0$), as we may as well assume for notational
convenience that the two sections differ in a neighborhood of 
$S_0$: here we get a {\it family} of surfaces where the sections
$\psi_{\cal F}$ and $\Phi$ disagree rather than
just a pair of points because the sections 
$\psi_{\cal F}$ and $\Phi$ are continuous. 

To set notation,
we rephrase this as follows: there is a family
$\{S_t\}$ of distinct Riemann surfaces and corresponding 
$\Gamma$-equivariant 
harmonic maps $f_t:\wt S_t\to T$,  Hopf differentials $\wt\Phi_t$, 
vertical foliations ${\cal F}_t$, and projections
$\pi_t:\wt S_t\to R_t$ to $\R$-trees with small $\Gamma$-actions and 
universal covering maps $p_t:\hyp\to S_t$
(choosing the notation so 
that $t=0$ corresponds to the original action).  Note that the trees
$R_t$ and morphisms $\phi_t$ are distinct, and that the 
foliations ${\cal F}_t$ represent different points in
$\PMF(S)$.  If this were not true then 
$\psi_{\cal F}(S_{t_1}) = \Phi(S_{t_1})$ for some $t_1>0$, which 
would force the sections $\psi_{\cal F}$ and $\Phi$ to agree over
$S_{t_1}$, contrary to the definition of the family $S_t$.

The heart of our argument is the case when the foliations 
${\cal F}_t$ are orientable and minimal.  We begin 
with a reduction towards that case.

\subsubsection{Some leaf is not closed}
\label{subsection:closedleaves}

Let $e\in E\subset T$ denote a point of $T$ which 
is not the image of a vertex in $R_0$ and which lies on 
the folded edge $E$ of $T$.  We consider the leaves of ${\cal F}_0$ 
containing 
$p_0\circ f^{-1}_0(e)\subset S_0$.  Since $e$ lies on the folded edge $E$ 
there are at least two of these.  Each such 
leaf which is a closed curve represents a (conjugacy class of) element of 
$\Gamma$ which fixes the edge $E\subset T$.

If each of these two leaves were closed, then they must represent the 
same element of $\pi_1(S)$: being simple closed curves, they do not 
represent powers of a common element of $\pi_1(S)$, hence some powers of 
these two elements in $\pi_1(S)$ must generate 
a free group since $S$ is closed and hyperbolic; this free 
subgroup of $\pi_1(S)$ stabilizes $E$, contradicting smallness.  But these 
two closed leaves are not even freely homotopic.  If they were then they 
would bound an annulus $A$ on $S_0$.  Since $A$ has 
Euler-characteristic zero and the boundary components are leaves of 
${\cal F}_0$, no singularity of ${\cal F}_0$ 
lies in $A$.  Hence the foliation ${\cal F}_0$ on this annulus would be
by closed curves parallel to the boundary and the harmonic map
$\pi\bigm|_A$ restricted to this annulus would map to an interval,
with constant boundary values. This forces the map to be everywhere
constant, so that the Hopf differential vanishes on $A$, hence
everywhere, an absurdity.

\subsubsection{The model case}
\label{subsection:model}

So we may assume that one of the components $\ell$ of
$p_0\circ f^{-1}_0(e)$ is not closed. Then consider a small arc 
$\alpha \subset S$ transverse to $\ell$ and to ${\cal
F}_0$. As the leaf $\ell$ is not closed, it is dense in a subsurface which we
might as well take to contain $\alpha$ (after maybe reducing the size
of $\alpha$ -see \cite{St}, \S11). Indeed, we can find a finite
number of edge points $e_1,\dots,e_n$ so that the trajectories
$p_0\circ f^{-1}_0(e_i)$ have closure equal to all of $S_0$.

Again, let $\alpha$ denote a small half-open arc transverse
to ${\cal F}_t$ on $S_t$ with its endpoint on the singularity
$q_0\in S$; we also assume that $f_t(\alpha)\subset E$, the
folded
edge, and that $\alpha$ is chosen small enough to ensure that
$\phi_t\bigm|_{\pi_t(\alpha)}$ is an isometry. 
By \ref{subsection:closedleaves}, we may 
assume that the nonsingular leaves through
$\alpha$ are not closed on $S_t$. 
(If a non-singular leaf were closed, it would be contained in a 
neighborhood of non-singular closed leaves (\cite{St}, \S9.3)
and so there would be no leaf through $\alpha$ which would also 
be dense in a subsurface containing $\alpha$. On the other hand,
if every neighborhood of $q_0$ in $\alpha$ had regular closed leaves,
since there are but a finite number of (maximal) ring domains
(i.e. maximal neighborhoods of regular closed leaves) in ${\cal F}_t$,
we see that a neighborhood of $q_0$ in $\alpha$ is contained in one of
these ring domains.  If this were true for all arcs $\alpha$ as above
with $f_t(\alpha) \subset E$, we would be in the situation of 
\ref{subsection:closedleaves}, a contradiction.)

We begin with the model case of ${\cal F}_t$ being orientable and
{\em minimal}, i.e., every non-singular leaf is dense.  The general case
will follow from technical modifications to the proof in this case,
but the essential ideas will be the same as in this model case. 

Now,
under the assumption that ${\cal F}_t$ is minimal and orientable, we see
that the first return map $P_t:\alpha\to\alpha$ determines an interval
exchange map $\sigma_t:\alpha\to\alpha$ on $\alpha$ (see \cite{St},
p.~58).  Moreover,
one can reconstruct the measured foliation $({\cal F}_t,\mu_t)$ directly
from the interval exchange map $\sigma_t:\alpha\to\alpha$. We recall that 
this
interval exchange map $\sigma_t$ is determined by looking at the largest
open subintervals $R_i(t)$ of $\alpha$ on which $P_t$ is continuous.  
The endpoints $\{x_0(t)=q_0,x_1(t),\dots,x_N(t)\}$ of these
subintervals are contained in singular leaves of ${\cal F}_t$, and hence
(have lifts to $\wt S$ which) project to vertex points of the tree $R_t$. 

We know that the set of vertex points in $R_t$ is totally
disconnected, as they are the image of the countable discrete set in 
$\hyp$ of zeroes of $\Phi_t$.  It is also easy to see from this that 
the set of vertex points of the tree $\phi_t(R_t)$ in $T$ is also 
totally disconnected.  
We now assume, postponing the proof until the end of this
subsection, that for each $t$ there is some vertex point 
$v \in R_t$ such that $\phi_t(v)$ a vertex point.  

\bigskip
\noindent
{\bf Continuity argument: }
Our main observation is that, since the $\Gamma$-equivariant maps
$f_t:\hyp^2\to T$ are continuous in $t$, we see that if 
$f_t(\wt{x_i(t)})$ is a vertex
in $T$, then as the vertices in $T$ are a totally disconnected set, the 
family
$f_t(\wt{x_i(t)})$ is constant in $t$. By the previous paragraph, there 
must 
be at least one endpoint $x_i(t)$
whose lift $\wt{x_i(t)}$ projects to a vertex in $R_t$.  
Since $\SF_t$ is minimal and $f_t$ is 
equivariant, we have that 
$\Gamma f_t(\wt{x_i(t)})=\Gamma f(\wt{x_i})$ is dense in $f(\wt{\alpha})$,
for lifts $ \wt{x_i(t)}$ and $\wt {\alpha}$ 
with $\wt{x_i(t)} \in \wt {\alpha}$. 
Letting 
$\Gamma_{x_i(t)} = \wt{\alpha} \cap
\pi_t^{-1}(\Gamma\pi_t\wt{x_i(t)})$, 
we see that $f_t\bigm|_{\Gamma_{x_i(t)}}$
is constant in $t$, which forces $f_t(\wt{x_j(t)})$ to 
be constant in $t$ for each $j$. 

Since the measure of
$\wt{\alpha}$ between consecutive vertices $\wt{x_i(t)}$ and 
$\wt{x_{i+1}(t)}$ (for
$i=0,\dots,N-1$) is determined by the distance $d_T$
($f_t(\wt{x_i(t)})$, $f_t(\wt{x_{i+1}(t)}$) in the tree $T$, we see that
these measures are also constant. Of course, after projecting
from the cover $\wt S$ to the surface $S$, we see that the endpoints
$x_i(t) \subset \alpha$ are also constant in $t$.

Finally, observe that the first
return maps $P_t:\alpha\to\alpha$ vary continuously in $t$ on the 
interiors of
the intervals $R_i(t)$
(and are affine there); hence, since the endpoints $x_i(t)$ are
constant in $t$, we see that the maps $P_t$ are constant in $t$ as well. We
conclude that the interval exchange maps $\sigma_t$ are constant in $t$,
so that $({\cal F}_t,\mu_t)=({\cal F}_0,\mu_0)$ after we reconstruct
$({\cal F}_t,\mu_t)$ from $\sigma_t:\alpha\to\alpha$. Hence we are done by 
Lemma \ref{lemma:key}.

\bigskip
\noindent
{\bf Proof that some vertex point maps to a vertex point: } 
Since this property is preserved under perturbations of the map, it is 
enough to prove the statement for some $t$.  

Suppose this were not the case. Then every vertex point
of every $R_t$ maps to an edge point of $T$.  Hence by 
Lemma \ref{lemma:benson} some edge 
point of each $R_t$ maps to a vertex point of $T$. Since there are 
finitely-many $\Gamma$-orbits of vertex points, there exists $\delta_t>0$ 
so that, on a $\Gamma-$fundamental domain of $R_t$,
any such edge point of $R_t$ has distance at least $\delta_t$ from 
any vertex point of $R_t$. For $t$ small, we may take all
$\delta_t > \delta$, for some fixed $\delta >0$.

We first claim that 
by making a small perturbation in Teichm\"{u}ller space
from $S$ to $S_t$ we may assume that ${\cal F}_t$ 
has a closed leaf $\lambda$
representing an edge point $x\in R_t$ within a $\delta/6$-neighborhood
of some vertex point $v_t$; necessarily, then 
there is a whole nondegenerate edge $E$ containing
$x$ which is both within a distance $\delta/3$ of 
$v_t$ and fixed by an element $g \in \Gamma$. 
This first claim
follows from essentially the same argument we used
in the continuity argument above: take a small arc 
which abuts the vertex $v_t \in R_t$, and consider the
image $\alpha$ on $S_t$ of a lift of that arc.  The foliation 
${\cal F}_t$ is determined by the interval exchange
defined by the first return map 
on that arc $\alpha$.  In particular, perturbations of 
${\cal F}_0$ are given by perturbations of that
first return map, and we can find such a perturbation ${\cal F}_0$
so that ${\cal F}_t$ has a closed leaf through $\alpha$.

Now we make a few observations about our situation:  since
(1) all vertices are being folded away to edge points creating
edges of radius at least $\delta/2$ from the image of 
$v_t$, but (2) on the surfaces $S_t$, no pair of adjacent sectors are 
having their $R_t$ images folded together (by the argument
late in the proof of Proposition 3.2), we see that for 
any point $e'$ in any edge $E'$
within $\delta/2$ of 
the image of $v_t$ in $T$, we must have at least
two distinct leaves on $S_t$ whose lifts project to $e'$.  But this 
contradicts smallness, as we showed in 
\S\ref{subsection:closedleaves}.  Hence some vertex point maps to a vertex 
point.
\bigskip

Next we begin to loosen the hypotheses of the model case so as to
eventually find ourselves in the general case, where ${\cal F}_t$ may be
non-orientable and have several minimal components. 

\subsubsection{Nonorientable case}
\label{subsection:general}
Let us first
drop the assumption that ${\cal F}_t$ should be orientable.  This is
merely a matter of generalizing the correspondence between measured
foliations $({\cal F}_t,\mu_t)$ and interval exchange maps $S_t$. The
idea here goes back to Strebel (see \cite{St}).  
We regard one side of $\alpha$ as
$\alpha^+$ and the other side as 
$\alpha^-$: if ${\cal F}_t$ is orientable, then
the rectangles $R_i(t)$ have one edge on $\alpha^+$ and another on
$\alpha^-$, but if ${\cal F}_t$ is not orientable, a rectangle may have
both edges on, say, $\alpha^+$. Yet, if we now regard the first return
map $P_t$ as a map $P_t:\alpha^+\cup\alpha^-\to\alpha^+\cup\alpha^-$, we 
can
consider an associated interval exchange map
$S_t:\alpha^+\cup\alpha^-\to\alpha^+\cup\alpha^-$ from which we can 
reconstruct
$({\cal F}_t,\mu_t)$. The endpoints $\{x_i(t)\}$ of the intervals
$R_i(t)$ on $\alpha^+\cup\alpha^-$ still 
(have lifts which) map continuously into the
disconnected set of vertices (constant in $t$) of $T$, so then, as
before the endpoints $\{x_i(t)\}$, and the map $P_t$, $S_t$ are
constant in $t$. We conclude that the measured foliations are also
constant in $t$.

\subsubsection{Breaking the model case into pieces}

We come finally to the most general part, where we no longer
require that ${\cal F}_t$ is minimal. Then for ${\cal F}_0$ choose a
collection of closed arcs $\alpha_1,\dots\alpha_n$ which are transverse to
${\cal F}_0$, and whose ${\cal F}_0$-orbits both cover $\hyp^2/\Gamma_0$ 
and
intersect at most along some compact singular leaves.
At this point, we also require the intervals $\alpha_i$ to have
corresponding interval exchange maps for ${\cal F}_0$ which are either
{\em irreducible}, i.e. we cannot (non-trivially) decompose
$\alpha_i=\alpha'_i\cup\alpha''_i$ with the interval exchange 
map $\sigma_i$ for $\alpha_i$
having a restriction $\sigma_i\bigm|_{\alpha'_i}:\alpha'_i\to\alpha'_i$ 
which
preserves the proper subinterval $\alpha'_i$, or correspond to a single
cylinder in ${\cal F}_0$ , so that the interval exchange is the identity
on a single cylinder. 

We claim that the measured foliation ${\cal F}_t$ on the
whole surface $\Sigma(t)$ is constant in $t$. This will give a
contradiction by Lemma \ref{lemma:key}, proving the theorem.

Let $\Sigma_i(t)$ be the subsurface of $S_t$ obtained by taking the
closure of the orbit of $\alpha_i$ along the leaves of ${\cal F}_t$. Our
restrictions on $\{\alpha_i\}$ have the effect of forcing either 
$\Sigma_i(0)$
to be a cylinder or a surface on which $\SF_0$ is minimal. 
We observe that
the argument given earlier for the cases where $\SF_0$ was minimal on the
closed surface $\hyp^2/\Gamma_0$ continue to hold for the case where 
$\SF_0$ is
minimal on $\Sigma_i(0)$. In particular, for $\Sigma_i(0)$ a subsurface 
with
almost every leaf dense, we see that the interval exchange maps 
$\sigma_i(t)$
must be constant in $t$. Yet, it is part of the basic construction of 
measured
foliations from interval exchange maps that the topology of $\Sigma_i(t)$ 
(as
well as $\SF_t$) is determined from the map $\sigma_i(t)$ (see, e.g.,
\cite{Ma1}). Thus, as $\sigma_i(t)=\sigma_i(0)$,
we see that $\Sigma_i(t)$ is homeomorphic to $\Sigma_i(0)$.

Now each boundary circle of each $\Sigma_i(t)$ is a
leaf of the foliation on that subsurface.  This leaf may be 
taken to be singular as it would otherwise be an interior leaf of a 
cylinder of
non-singular homotopic leaves, counter to the construction of 
$\{\alpha_i\}$. Thus the continuity argument also shows that the 
foliations 
on the cylindrical subsurfaces $\Sigma_j(t)$ are constant in $t$.  
Hence the measured foliation on each subsurface $\Sigma_i(t)$ is
constant in $t$. 
Finally, whenever two subsurfaces
$\Sigma_i(t)$ and $\Sigma_j(t)$ have a common boundary component $C(t)$, 
the
continuity argument shows that $C(t)$ cannot become a cylinder at any time 
$t$
as this would require the single vertex $f_t(C(t))$ to continuously deform 
into
a non-trivial family of pairs of vertices, an absurdity. So we see that the
identification of all the boundary components of all the
$\Sigma_i(t)$ are constant over $t$, so that
${\cal F}_t$ is constant.

\section{Appendix}

This appendix is dedicated to a proof of Proposition 
\ref{proposition:moving}, which is partly implicit in \cite{HM} and partly 
a ``folklore theorem''.  We provide here an elementary, geometric, and 
self-contained proof due almost entirely
to Howard Masur (personal communication), who graciously permitted
us to reproduce it here.

The proof can be reduced to the following claim: If either 

\begin{enumerate}
\item  $q$ has a pair $\{z_1,z_2\}$ of distinct zeros 
connected by an arc $A$ of a leaf of its vertical foliation, or
\item $q$ has a $k$-pronged singularity at $z_3$, 
and $2$ arbitrary sectors $s_1,s_2$ of this $k$-prong are specified, 
\end{enumerate}

then there is a Riemann surface $\SR^*$ and a
holomorphic quadratic differential $q^*$ on $\SR^*$ so that the vertical
foliation of $q^*$ is measure equivalent to $(\SF,\mu)$ and
\begin{enumerate}
\item (in case (1) above) the zeros of $q^*$
corresponding to $\{z_1,z_2\}$ coincide, or 
\item (in case (2) above) the images of the sectors $s_1,s_2$ under the
  equivalence of foliations meet along an arc.
\end{enumerate}

The proposition follows from the claim as follows.  First apply (1)
above to get $s_1$ and $s_2$ as sectors abutting a common 
singularity $z$.  Then apply (2) above and we are done.

We are left to prove the claim.  

\bigskip
\noindent
{\bf Single cylinder case. } We first prove the claim for {\em Jenkins
  differentials}, i.e. those differentials whose vertical 
foliations are but one foliated open 
(right Euclidean) cylinder $C$ with all singularities lying 
in $\partial \overline{C}$.  Here $S$ can be thought of as an
identification space $\pi:\overline{C}\rightarrow S$, with
identifications being made on $\partial \overline{C}$.  Let $C_1,C_2$
denote the 2 components of $\partial \overline{C}$.  Note that the graph 
$L$ lies in $\partial \overline{C}$, and all the 
singularities on a single component of $\partial \overline{C}$ are 
connected by 
$L$. Moreover, there are natural correspondences between
topological or geometric operations on the surface $S$ and 
topological or geometric operations on $\overline{C}$. This means
that if we continuously deform $C$ to another right Euclidean 
cylinder $C^*$ (so that  there is a canonical correspondence
of identifications on  $\partial \overline{C^*}$), then the 
canonical quadratic differential $q^*$ on $C^*$ (defined so that
the metric $|q^*|$ agrees with the Euclidean metric on $C^*$ and all
of whose vertical leaves are parallel to  
$\partial \overline{C^*}$) descends
to a quadratic differential $q^*$ on the identified surface $S^*$
with the vertical foliation of $q^*$ on $S^*$ being Whitehead equivalent
to the vertical foliation of $q$ on $S$.

To prove (1) and (2) above, we will first perform the desired
operation on $C$ to obtain a new Euclidean 
cylinder $C^\ast$ with canonically
determined quadratic differential $q^*$ as above.  
The important thing to check in each case is that we can do this so 
that the resulting Euclidean 
lengths $\ell(C^\ast_1)$ and $\ell(C^\ast_2)$ 
of the 2 components of $\pov{C^\ast}$ 
are equal.  This imediately implies that the identification $\pi$
determines an identification $\pi^\ast:C^\ast\rightarrow S^\ast$ to a
Riemann surface $S^\ast$, and that the canonical quadratic differential on
$C^\ast$ descends to a quadratic differential $q^\ast$ on $S^\ast$.
By construction $q^\ast$ has 
vertical foliation measure equivalent to that of the
vertical foliation of $q$.

Consider case (1).  Let $A_1,A_2\subset \partial \overline{C}$ 
denote the 2 components of $\pi^{-1}(A)$, where we recall that
$A$ is the arc of the vertical foliation we wish to collapse.
Note that $\ell(A_1)=\ell(A_2)$.

\medskip
\noindent
{\bf Case 1a: $A_1$ and $A_2$ lie in different components of $\partial
  \overline{C}$.} In this case contract both $A_1$ and $A_2$ to a
point to give a Euclidean cylinder $C^\ast$.  Since $C^\ast_1$ and
$C^\ast_2$ have the same Euclidean length, so we are done by the
above.  

\medskip
\noindent
{\bf Case 1b: $A_1,A_2$ lie on the same component of $\pov{C}$. }
First note that the arcs of $L$ have 
preimages in $\partial \overline{C}$ which come in 
pairs, as neighborhoods of the arcs on the identification
space have full neighborhoods, while neighborhoods of
arcs on $\partial \overline{C}$ have only half-neighborhoods.
Hence since $A_1,A_2$ lie on the same component of $\pov{C}$, 
we must be able to find some collection of pairs of
arcs on the other component the sum of whose lengths is
at least that of the sum of the lengths of 
$A_1$ and $A_2$. (This is just a pigeon-hole principle:
the arcs come in pairs whose lengths are equal and for which
total lengths of all the arcs 
are the sum of the lengths of the boundary
components of $\partial \overline{C}$, yet each of 
these boundary components have the same lengths, so the fact 
that $A_1$ and $A_2$ contribute solely to one component of
$\partial \overline{C}$ forces some other family of
pairs to contribute at least as much solely to the
other component of $\partial \overline{C}$.)
Thus we act as before, contracting $A_1$ and $A_2$ on
one component of $\partial \overline{C}$ and simultaneously
some other pairs of arcs the same amount on the
other component of $\partial \overline{C}$. It is
quite important here that the contraction of the
other components has no effect on our claim or our
purpose; the proof of the first part of the claim
concludes as before.
\bigskip

Now to prove part (2) of the claim.  Under the identification map
$\pi:C\rightarrow S$, each sector $s_i, i=1,2$ has a
unique pre-image on $C$ as a neighborhood $U_i$ of a vertex $v_i$.

\medskip
\noindent
{\bf Case 2a: $v_1$ and $v_2$ lie on different components of $\pov{C}$.}
 
We split the
vertex $v_1$ into a pair of vertices $v_{1,1}$ and $v_{1,2}$
connected by an arc $A_1$, and we
split the vertex $v_2$ into a pair of vertices $v_{2,1}$ and $v_{2,2}$
connected by an arc $A_2$ of the same length as $A_1$. We then 
re-identify the cylinder as before, with the only changes being that
instead of identifying $v_1$ to $v_2$, we send $A_1$ isometrically
onto $A_2$
(there is a unique way to do this which preserves the ordering of the
sectors). The resulting surface gives $S^*$ and $q^\ast$ as required. 

\medskip
\noindent
{\bf Case 2b: $v_1$ and $v_2$ lie on the same component, say $C_1$, of
  $\pov{C}$. } Split $v_1$ and $v_2$ as in Case 2b.  We now do a 
further deformation to make $\ell(C^\ast_1)=\ell(C^\ast_2)$.  

If for some compact
singular arc $B \subset L$, we have
both components of $\pi^{-1}(B)$  lying
on $C_2$, then by lengthening $B$ we could
achieve $\ell(C^\ast_1)=\ell(C^\ast_2)$.  If this isn't true, then 
by the pigeon-hole priciple, for each such $B$ we know  
$\pi^{-1}(B)$ has one component on $C_1$ and one on $C_2$.

Now observe that any singularity 
on the surface $\SR$ with, say, $k$ sectors, admits a cyclic
ordering of these sectors $s_1, ...s_k$ (where the closure
of $s_2$ meets the closure of $s_1$ on one side and the closure of 
$s_3$ on the other side, and so on).  Since we are in a case where
each edge incident to a singularity on $\SR$ has preimages on both
boundary components $C_1$ and $C_2$ of $\pov{C}$,
and since sectors have preimages near components of 
$\partial \overline{C}$ where their bounding arcs have preimages, 
we see that the sectors $s_1,\ldots ,s_k$ 
also alternate between having preimages in $C_1$ and in $C_2$. 
This implies that all of the singularities on the surface  $\SR$ have 
an even number of sectors.

We now claim that there are 
vertices $w_1,w_2$ in $C_2$ which are still identified
by the identification rules, even after the splitting of 
$v_1$ and $v_2$.  
(Here the subtlety is that by first splitting $v_1$ and $v_2$, 
we have changed the identification rules, and hence the 
orbits of identified vertices on  $\partial \overline{C}$. Our
vertices $w_1$ and $w_2$ must not only then correspond to each
other by the original identification rules, but they must also
lie in the same new orbit of vertices on  $\partial \overline{C}$,
after the splitting of $v_1$ and $v_2$.)
This finishes
the proof of case 2b since we then split $w_1$ and $w_2$ to make
$\ell(C^\ast_1)=\ell(C^\ast_2)$.  

To see that there are such vertices $w_1$ and $w_2$, we 
recall that the total multiplicity
of zeroes of a holomorphic quadratic differential on a Riemann
surface $\SR$ of genus $g$ is equal to $4g-4$ (Riemann-Roch). 
Thus, since in the 
case under consideration all of the singularities have an even number
of sectors (and hence an even order of zero), we see that there is
either one singularity $z_0$ with at least six sectors, or several
singularities which all have at least four sectors.  In the first
case we see that any initial splitting of $z_0$ 
(by splitting a pair of vertices $v_1$ and $v_2$ on
the same component $C_1$ of $\partial \overline{C}$) would leave a
topological foliation with two singularities of which
at least one would have four sectors, with two of those
sectors
having preimages on $C_2$: we would then split 
the vertices of those sectors, say $w_1$ and $w_2$
to finish the case. In the second case, there is at the outset 
a singularity on $\SR$ whose preimages do not include $v_1$ and $v_2$,
and which has at least a pair of sectors with preimages on
$C_2$, as desired.

\bigskip
\noindent
{\bf General case. }
We prove the general case by the now standard technique of approximating.  
By \cite{Ma2} we may approximate $q$ by a 
 sequence $\{q_n\}$ of Jenkins differentials on $\SR$. In case (1), 
 let $A$ denote the
arc of the vertical foliation of $q$ which we wish to contract.
As $q_n$ approximates $q$, there is an arc $A_n\subset\p\ov C_n$ which
approximates $A$. Furthermore, there is a contractible 
neighborhood $U$ on the underlying
differentiable surface which is a neighborhood of the arc $A$ and all arcs
$A_n$ for $n$ sufficiently large.

Now, the Riemann surface $\SR$ is an 
identification space of each cylinder $\ov
C_n$, with identifications being made on $\p\ov C_n$. As in the 
single cylinder case above, we can form new Riemann
surfaces $\SR^*_n$ equipped with 
quadratic differentials $q^*_n$ by contracting
the arcs $A_n\subset\p\ov C_n$ and 
identifying as before; here the arc $A_n$ on
$\SR_n$ bounded by a pair of low 
order zeros is replaced on $\SR^*_n$ by a single
higher order zero, say $z^*_n$. 

The important thing to notice about
this operation is that the complement $V=U^c$ of the neighborhood
$U\subset\SR$ is approximated by 
the closure of an open
set $V_n$ on $\SR^*_n$ which only avoids a
small neighborhood of the high order zero $z^*_n$;
moreover, the conformal structures on
$V_n$ compare uniformly to the conformal structure on $V$, hence to each 
other.  Hence, by passing to a subsequence if necessary, we have that 
that $\SR^*_n$ converges in (the interior of) Teichmuller space 
to a Riemann surface $\SR^*$. It also
follows that $q^*_n$ converges to a holomorphic quadratic differential, say
$q^*$, on $S^*$; as $q_n$ approximates $q$ and $q^*_n$ is measure equivalent to
$q_n$, we see that $q^*$ is measure equivalent to $q$. Moreover, as the 
foliation of $q^*_n$ results from a Whitehead move 
applied to the foliation of $q_n$
(which contracts $A_n$ to a point), the foliation of $q^*$ is obtained from the
foliation of $q$ via a Whitehead move which contracts $A$.  

Case (2) 
is virtually identical: we still have uniform convergence of the
conformal structures outside the pair(s) of neighborhoods of the vertices
(or arcs) we are splitting to pairs of vertices connected by an arc.

\endproof

\noindent
Benson Farb:\\
Dept. of Mathematics, University of Chicago\\
5734 University Ave.\\
Chicago, Il 60637\\
E-mail: farb@math.uchicago.edu
\medskip

\noindent
Michael Wolf:\\
Dept. of Mathematics, Rice University\\
Houston, TX 77251\\
E-mail: mwolf@math.rice.edu

\end{document}